\documentclass{article}

\usepackage{theorem}
\usepackage{amssymb}

\def\sf{{\cal F}}

\def\R{{\mathbb R}}
\def\N{{\mathbb N}}
\def\C{{\mathbb C}}
\def\rd{{\rm d}}
\newcommand{\QED}{\hfill$\Box$\par}
\def\sskip{\vskip .2cm}
\def\jtext#1{{~\textrm{#1}~}}

\newtheorem{theorem}{Theorem}[section]

\title{Convergence from below suffices}
\author{J.~F.~Feinstein}
\date{}
\begin{document}

\maketitle
\begin{abstract}An elementary application of Fatou's lemma gives a strengthened version
of the monotone convergence theorem. We call this the {\em convergence from below theorem}.
We make the case that this result should be better known, and deserves a place in
any introductory course on measure and integration.
\end{abstract}
\section{The convergence from below theorem}
Three famous convergence-related results appear in most introductory courses on measure and integration:
the monotone convergence theorem, Fatou's lemma and the dominated convergence theorem.
In teaching this material it is common to follow the approach taken in, for example, \cite[Chapter 1]{Rudin}.
There Rudin begins by proving the monotone convergence theorem and then deduces Fatou's lemma. Finally, he deduces the dominated convergence theorem from Fatou's lemma. The result which we call the {\em convergence from below theorem}
(Theorem \ref{CBT} below) is essentially distilled from this proof of the dominated convergence theorem (\cite[pp. 26-27]{Rudin}).
We do not claim originality for this result, or for the related Theorem \ref{infinite-integral}. They are presumably known, although we know of no explicit references for them. However, we wish to make a case that that they should be better known than they are. In particular, we suggest that Theorem \ref{CBT} deserves a name and a place in the syllabus when this material is taught.
\sskip
Throughout we discuss results concerning pointwise convergence. In the usual way, there are versions of all these results in terms of almost-everywhere convergence instead.

For convenience, we shall use the following terminology.
Let $X$ be a set, let $(f_n)$ be a sequence of functions from $X$ to $[0,\infty]$ and let
$f$ be another function from $X$ to $[0,\infty]$. We say that the functions $f_n$ converge to $f$
{\bf from below on $X$} if the functions $f_n$ tend to $f$ pointwise on $X$ and
$f_n(x) \leq f(x)$ $(n \in \N,~x \in X)$.
We say that the functions $f_n$ converge to $f$
{\bf monotonely from below on $X$} if the functions $f_n$ tend to $f$ pointwise on $X$ and, for all $x \in X$, we have
$f_1(x) \leq f_2(x) \leq f_3(x) \leq \cdots$.

\medskip

We begin by recalling the statement of the monotone convergence theorem.
\newpage
\begin{theorem}
{\bf (Monotone convergence theorem)}
Let $(X,\sf,\mu )$
be a measure space, and let $f : X \to [0,\infty ]$ be a measurable function.
Let $(f_n)$ be a sequence of measurable functions from $X$ to $[0,\infty]$ which converge to $f$ monotonely from below on $X$.
Then
\[
\int _X f\,\rd \mu = \lim _{n\to \infty} \int _X f_n\,\rd \mu\,.
\]
\end{theorem}

The measurability assumption on
$f$ is, of course, redundant here as it follows from the pointwise convergence of $f_n$ to $f$.
We now observe that an elementary application of Fatou's lemma shows that we may weaken the monotone convergence assumption. We have not found this result stated explicitly in the literature, and it does not appear to have a name. We propose to call it the {\em convergence from below theorem}.

The concepts involved in the statements and applications of the monotone convergence theorem and the dominated convergence theorem are relatively simple. We suggest that convergence from below is a similarly simple concept, which should appeal to all levels of student. In particular, those students who find the concepts of $\liminf$ and $\limsup$ difficult may be happier applying the convergence from below theorem rather than Fatou's lemma (where possible).
\begin{theorem}
{\bf (Convergence from below theorem)}
\label{CBT}
Let $(X,\sf,\mu )$
be a measure space, and let $f : X \to [0,\infty ]$ be a measurable function.
Let $(f_n)$ be a sequence of measurable functions from $X$ to $[0,\infty]$ which converge to $f$ from below on $X$.
Then
\[
\int _X f\,\rd \mu = \lim _{n\to \infty} \int _X f_n\,\rd \mu\,.
\]
\end{theorem}
{\bf Proof.}
Clearly
\[
\limsup_{n\to \infty}\int_X f_n\,\rd \mu \leq \int_X f\, \rd \mu\,.
\]
However, by Fatou's lemma,
\[
\int_X f\,\rd \mu \leq \liminf_{n\to\infty}\int_X f_n\,\rd \mu\,.
\]
The result follows immediately.
\QED
\sskip
\noindent
{\bf Remarks.}
\begin{enumerate}
 \item[(1)]
The monotone convergence theorem is now a special case of our stronger convergence from below theorem.
\item[(2)]
In the case where $\int_X f\,\rd \mu < \infty$, the convergence from below theorem is an immediate consequence of the dominated convergence theorem.
\item[(3)]
In the case where $\int_X f\,\rd \mu = \infty$, the result does not follow directly from either the monotone convergence theorem or the dominated convergence theorem. The following elementary result clarifies the situation in this case.
\end{enumerate}

\begin{theorem}
\label{infinite-integral}
Let $(X,\sf,\mu )$
be a measure space, and let $f : X \to [0,\infty ]$ be a measurable function with $\int_X f\,\rd \mu = \infty$.
Let $(f_n)$ be a sequence of measurable functions from $X$ to $[0,\infty]$ which converge to $f$ pointwise on $X$.
Then
\[
\lim _{n\to \infty} \int _X f_n\,\rd \mu\, = \infty.
\]
\end{theorem}
{\bf Proof.}
By Fatou's lemma,
\[
\infty = \int_X f\,\rd \mu \leq \liminf_{n\to\infty}\int_X f_n\,\rd \mu\,.
\]
It follows immediately that
$\lim _{n\to \infty} \int _X f_n\,\rd \mu\, = \infty$, as required.
\QED
\sskip
We suggest that the convergence from below theorem deserves a place between Fatou's lemma and the dominated convergence theorem: the dominated convergence theorem may be deduced from the convergence from below theorem as follows. This proof is based on the proof given in \cite[pp. 26-27]{Rudin}, but applying the convergence from below theorem in the middle.

\begin{theorem} {\bf (Dominated convergence theorem)}
Let $(X,\sf,\mu )$
be a measure space, let $g : X \to [0,\infty]$ be a measurable function.
with $\int_X f\,\rd \mu < \infty$ and let $f$ be a measurable function from $X$ to $\C$.
Let $(f_n)$ be a sequence of measurable functions from $X$ to $C$ which converge to $f$ pointwise on $X$
and such that $|f_n(x)| \leq g(x)$ $(n \in \N, x \in X)$.
Then
\[
\lim_{n \to \infty}\int_X|f_n-f|\,\rd \mu\, = 0
\]
and
\[
\int _X f\,\rd \mu = \lim _{n\to \infty} \int _X f_n\,\rd \mu\,.
\]
\end{theorem}
{\bf Proof.}
The second equality follows quickly from the first.
To prove the first equality, observe that the non-negative, measurable functions $2g-|f_n-f|$ converge to the function $2g$ from below. Thus, by the convergence from below theorem,
\[
\lim_{n \to \infty}\int_X\left( 2g -|f_n-f|\right)\,\rd \mu\, = \int_X 2g \,\rd \mu\,.
\]
The result now follows by subtracting $\int_X 2g \,\rd \mu\,$ from both sides and rearranging.
\QED
\sskip

As discussed above, the convergence from below theorem is more than covered by a combination of the dominated convergence theorem and Theorem \ref{infinite-integral}. Also, since the convergence from below theorem is such an elementary consequence of Fatou's lemma, any applications may also be deduced from that lemma.
However, the monotone convergence theorem continues to be used in the literature, and any application of the monotone convergence theorem can be replaced directly by an application of the
convergence from below theorem.
Of course, we then only need to check the weaker conditions of the latter theorem.

Also, the convergence from below theorem can be used to give elegant solutions to simple problems where neither the monotone convergence theorem nor the dominated convergence theorem apply directly. Here is such an application (an elementary undergraduate exercise).
\sskip
\noindent
{\bf Exercise.}
Let $\lambda$ denote Lebesgue measure on $\R$.
Prove that, for every Lebesgue measurable subset $E$ of $\R$, we have
\[
\int_E x^2\, \rd \lambda(x) = \lim_{n \to \infty} \int_E \left( x^2 - \frac{1}{n} | x \sin n x|\right) \, \rd \lambda(x)\,.
\]
\sskip
\noindent
{\bf Solution.}
Since $|x \sin n x| \leq  n x^2$ ($n \in \N$, $x \in \R$), the result is an immediate consequence of the convergence from
below theorem.

We may, instead, apply Fatou's lemma directly. This does, of course, lead to a quick solution which essentially proves the convergence from below theorem again along the way.

We may also consider separately the cases where $\int_E x^2\, \rd \lambda(x)<\infty$ and
where $\int_E x^2\, \rd \lambda(x)=\infty$.
In the first case we may apply the dominated convergence theorem, and in the second case we may use Theorem \ref{infinite-integral}. However the use of the convergence from below theorem renders this splitting into two cases unnecessary.

\section{Proving the convergence from below theorem directly}
Above we suggested following the usual development of the theory, but inserting the convergence from below theorem between Fatou's lemma and the dominated convergence theorem. There are several alternatives, however. For example, we can prove Fatou's lemma directly first and then deduce the convergence from below theorem. The monotone convergence theorem and the dominated convergence theorem then follow easily.

Another approach is to modify the standard proof of the monotone convergence theorem (\cite[1.26]{Rudin}) in order to give a direct proof of the convergence from below theorem. The monotone convergence theorem, dominated convergence theorem and Fatou's lemma are then corollaries of this.
We conclude with such a direct proof.

In this proof we
avoid explicit reference to $\liminf$ and $\limsup$ in order to make the proof more accessible to students who have
difficulty with these concepts. However, only minor changes are needed to give a direct proof of Fatou's lemma instead.
\sskip
\noindent{\bf Direct proof of Theorem \ref{CBT}.}
First note that we have $\int_X f_n\,\rd \mu \leq \int_X f\,\rd \mu$ $(n \in \N)$. Thus it is sufficient to prove that,
for all $\alpha < \int_X f\,\rd \mu$, $\int_X f_n\,\rd \mu$ is eventually greater than $\alpha$, i.e., there is an $N \in \N$ such that, for all $n \geq N$, we have $\int_X f_n\,\rd \mu > \alpha$.
Given such an $\alpha$, the definition of the integral tells us that there is a nonnegative, simple measurable function $s$ with $s(x)\leq f(x)$ $(x \in X)$ and such that
$\int_X s \,\rd \mu > \alpha$.
Choose $c \in (0,1)$ large enough that
$\int_X c s \,\rd \mu > \alpha$.
Set
$A_n = \{x \in X: cs(x) \leq f_n(x)\}$
and, for each $k \in \N$, set
\[
B_k = \bigcap_{n \geq k} A_n = \{x \in X: cs(x) \leq f_n(x) \jtext{for all} n \geq k\}.
\]
Clearly, $B_1\subseteq B_2\subseteq \cdots$. We claim that $\bigcup_{k=1}^\infty B_k = X$.
Let $x \in X$. If $s(x)>0$, then $cs(x)<f(x)$, and so
$x \in B_k$ provided that $k$ is large enough. On the other hand, if $s(x)=0$, then $x \in B_k$ for all $k \in \N$.
This proves our claim.
By standard continuity properties of measures, we have
\[
\int_X{cs\,\rd \mu} = \lim_{k\to \infty} \int_{B_k} cs\,\rd\mu\,.
\]
Choose $N \in \N$ such that $\int_{B_N} cs\,\rd\mu > \alpha$.
For all $n \geq N$ and $x \in B_N$ we have $cs(x) \leq f_n(x)$.
Thus, for $n \geq N$, we have
\[
\int_X f_n \,\rd\mu\, \geq \int_{B_N} f_n \,\rd\mu\, \geq  \int_{B_N} cs \,\rd\mu\, > \alpha,
\]
as required.
\QED

School of Mathematical Sciences

 University of Nottingham

 Nottingham NG7 2RD, UK

 email: Joel.Feinstein@nottingham.ac.uk

\end{document}